# Parametric transformation functions in the Kaprekar routine (I)

Fernando Nuez

Retired professor of Universidad Politécnica de Valencia (Spain)

**Abstract**

The Kaprekar transformation is uniquely determined by **α** parameters based on differences between symmetric values in the ordered numeric sequence. The parametric analysis of the iteration of the process requires developing $K_i$ functions that provide the parameters of the transformed number $K_i\ (\boldsymbol{\alpha}) = \boldsymbol{\alpha}$'. This is the objective of this work. As application examples we use these functions to deduce constants and cycles. But its main utility is to study the algebraic architecture of the transformation trees that we develop in the second part of this work.

1. **Introduction**

In 1949 the Indian mathematician Dattatreya Kaprekar discovered what is known as Kaprekar's constant, number 6174 (Kaprekar, 1949). He defined the following transformation:

Take any 4-digit number n, and sort its digits in a descending order, which returns another number X. X's digits are inverted, yielding Y. Then, n' = X-Y is calculated, and the process is repeated as many times as possible until it is unnecessary, as the result is always 6174. Thus, if n = 8082, X = 8820, Y = 0288, n' = X-Y = 8532, X' = 8532, Y' = 2358, n'' = 6174, X'' = 7641, Y'' = 1467, n''' = X''- Y'' = 6174, $n^{4)}$ = 6174, …When applying the rutine to the set of three-digit numbers, the result is 495.

Gardner (1975) contributed to popularize the figure of Kaprekar and the number 6174. In his mathematical games column in Scientific American he focuses on generated and self-generated numbers, but dedicates a paragraph to the Kaprekar constant.

Since then, much research has been done on various aspects of the process: number of digits, numerical bases, existence of constants and cycles, maximum distances,… A good bibliographic review can be seen in Yamagami's (2018). However, there are still many information gaps on algebraic aspects of the Kaprekar routine, and today there is renewed interest in the subject.

We are especially interested in the algebraic architecture of transformation trees and in existing symmetries. For this we have approached the process using as a tool the $K_i$ parametric transformation functions that we develop in this work. In the second part of the same, these functions are used to study symmetries in cases of 2, 3, 4 and 5 digits and a general methodology is presented (Nuez, 2021).

## 2. Formalization of the process for 4-digit and base 10

Numbers will be represented according to thousands, hundreds, tens and ones as usual in the decimal numeral system.

$n = a \times 10^3 + b \times 10^2 + c \times 10^1 + d \times 10^0 = a\,b\,c\,d = (a\,b\,c\,d)$

When necessary, a parenthesis will be used for clarity purposes. The reference set $A_4$ will include numbers up to four digits long, except for the zero and the nine numbers where all four digits are the same.

$A_4 = \{n \in \mathbb{Z}, 0 < n < 9999; n \neq 1111 \times s, s = 0, 1, \ldots 9\}$

Numbers with less than four digits will be represented with 4 digits by adding zeros to the left. Thus, $1 = 0\,0\,0\,1$. Similarly, we will later on analyze the sets $A_w$ of w digits, for w = 2, 3 and 5.

We will introduce operator O, which sorts digits in a descending order. Thus, $O\,(0\,0\,0\,1) = 1\,0\,0\,0$. Any permutation P of a given number's digits returns the same result (image) when O is operating. For example,

$n = 2\,5\,0\,9, \; m = P\,(n) = 5\,2\,9\,0, \; O(n) = O\,(m) = 9\,5\,2\,0$

The number having all its digits sorted in a descending order will be represented by $n_0 = O\,(n) = (a_0\,b_0\,c_0\,d_0)$

The parameters $\alpha$ and $\beta$ will be defined on the sequence of ordered digits. These parameters are essential to analyze the dynamics of Kaprekar's routine.

$O\,(n) = (a_0\,b_0\,c_0\,d_0), \quad a_0 \geq b_0 \geq c_0 \geq d_0, \quad a_0 > d_0, \quad \alpha = a_0 - d_0, \quad \beta = b_0 - c_0,$
$n_0 = (d_0 + \alpha \;\; c_0 + \beta \;\; c_0 \;\; d_0)$ [1]

Since $c_0 \geq d_0$ and $d_0 + \alpha \geq c_0 + \beta \rightarrow \boldsymbol{\alpha \geq \beta} \quad 0 < \alpha \leq 9, \; 0 \leq \beta \leq 9$ [2]

We will write $p\,(n) = (\alpha, \beta)$, which is to be read as "parameters of n: $\alpha$ and $\beta$"

The operation of sorting the digits in descending and ascending order and subtracting the second from the first will be represented by K. Naturally, operator O is a part of K.

$K\,(n) = K\,(a\,b\,c\,d) = (a'\,b'\,c'\,d') = n'$

The digits of n' are not necessarily sorted. Thus,

$n = 2\,5\,0\,9, \; n_0 = O\,(n) = 9\,5\,2\,0 \quad p\,(n) = (9-0, 5-2) = (9,3)$

$n' = 9\,5\,2\,0 - 0\,2\,5\,9 = 9\,2\,6\,1, \; p\,(n') = p\,[O\,(n')] = (9-1, 6-2) = (8,4)$

Insisting on the notation,

$n = (a\,b\,c\,d): a = 2, \; b = 5, \; c = 0, \; d = 9$

O (n) = (a$_0$ b$_0$ c$_0$ d$_0$): a$_0$ = 9,   b$_0$ = 5,   c$_0$= 2,   d$_0$ = 0

n' = (a' b' c' d'): a' = 9,   b' = 2,   c' = 6,   d' = 1

In the case n = 0 0 0 1,  p(n) = (1,0), n' = 0 9 9 9,  p (n') = (9,0)

In general

K(n) = (d$_0$+α   c$_0$+β   c$_0$   d$_0$) – (d$_0$   c$_0$   c$_0$+β   d$_0$+α) = (α   β   –β   –α)

As all digits must be natural numbers between 0 and 9, both included, negative digits will be converted to positive by moving units from a higher range (thousands, hundreds, tens) to a lower one (hundreds, tens, ones). Thus,

- If β > 0
  n'= K(n) = (α   β-1   9- β   10- α) = f$_1$ (α,ß) [3]
  a'= α,   b'= β-1,   c'= 9- β,   d'= 10- α
  S (n') = a' + b' + c' + d' =18 → n' = $\dot{9}$
  Thus,  n = 2 5 0 9,  p (n) = (9,3),  n' = (9   3-1   9-3   10-9) = 9261
- If β=0
  n'= K (n) = (α-1   9   9   10- α) = f$_2$ (α,0) [4]
  S (n') = 27 → n' = $\dot{9}$
  Thus, n = 0001,  p (n) = (1,0),  n' = (1-1   9   9   10-1) = 0999

These expressions are a particular case of those developed by Prichett et al. (1981). They formalize Kaprekar's routine and highlight the importance of parametrization.

  a) The way the digits appear is not important, due to the necessary descending and ascending arrangement before subtracting. Naturally, all permutations P$_i$ of the digits in a number, up to 4! = 24 different numbers if all the digits are different, will verify O (P$_i$) = O (P$_j$)  ∀ i,j and will return the same image K [P$_i$ (n)] = K (n), i=1, 2 ,… 24.
  b) If p (m) = p (n) and m ≠ P$_i$ (n) → K (m) = K (n)
     That is, regardless of what the digits of two numbers are, even non-permuted digits, if their parameters are the same, both will have the same image and vice versa.
     Thus, let us recall
     n = 2509,  p (n) = (9,3),  K (n) = 9261
     since m = 1904, O (m) = (9  4  1  0)  p (m) = (9,3)  = p (n)
     then K (m) = 9261 = K (n)
     In addition, all numbers which are permutations of n when converted will yield 9261. The same will be true for permutations of m. But m and n have different digits.
  c) If p (m) = (α$_1$,ß$_1$) and p (n) = (α$_2$,ß$_2$),   K (m) = K (n) ↔ α$_1$ = α$_2$,  ß$_1$ = ß$_2$ [5]
  d) If O (m) = O (n)  → K (m) = K (n). The reciprocal is not true.

3. **Absorbing set B$_4$**

   We will call B$_4$ ⊂ A$_4$ to the set of those numbers transformed by K
   B$_4$ = {n' = K (n); n ∈ A$_4$} [6]

By its own definition, this is a closed set: n'' = K (n') ∈ $B_4$. Its nature conditions the existence of certain symmetries.

Going back to expressions [3] and [4], the image n' in $B_4$ of any n, p(n) = (α, ß) must necessarily be a multiple of 9, n' = $\dot{9}$ , and satisfy the following conditions:

-If ß > 0

a' + d' = α + (10-α) = 10 ; b' + c' = (ß-1) + (9-ß)=8

-If ß = 0

a'+d' = 9 = b' = c'

These conditions are necessary but not sufficient for a number satisfying them to belong to $B_4$.

Thus, n = 4446 satisfies the first restriction a' + d' = 4+6 = 10, b'+c' = 8 but does not belong to $B_4$, as there is no m ∈ $A_4$; p (m) = (α,ß) such that K (m) = n. Indeed, K (m) = (α  ß-1  9-ß  10-α) = 4446 → α = 4, ß = 5 which is not possible, since by [1] and [2], α ≥ ß

Ultimately, n ∈ $A_4$, n ∉ $B_4$

Similarly, n = 9 9 9 0 satisfies the second restriction but n ∉ $B_4$ since ∄m, K(m) = (α-1  9  9  10-α) = 9 9 9 0 → α = 10, impossible since 0<α≤9

However,
- If n = (a  b  8-b  10-a)  0 ≤ a ≤ 9, a ≥ b+1 → n ∈ $B_4$
  since ∃ m ∈ $A_4$, p (m) = (a, b+1), K (m) = n
- If n = (a  9 9 9-a)  a ≤ 8 → n ∈ $B_4$
  since ∃ m ∈ $A_4$, p(m) = (a+1, 0), K (m) = n

Shortly, all, and only those numbers n= (a b c d) satisfying any of the two following requirements are images of K – and thus belong to the absorbing set $B_4$ –:

Requirement 1) a+d = 10, b+c = 8,   a ≥ b+1                                    [7]

Requirement 2) a+d = 9 = b = c,   a ≤ 8                                         [8]

since in both cases there is some m with β > 0 or β = 0  respectively, which satisfies K (m) = n.

4. **The parametric transformation functions $K_i$**

Let us use the following notation:
- K is used in a generic sense to refer to "image of". Its argument can be a number or the parameters of a number, yielding numerical images or parametrical images, respectively

Thus, K (n)=n';   K (α,ß) = (α',ß'),  p (n) = (α,ß),   p (n') = (α',ß')

- If functions are used specifically, all of them act upon parameters, but it is necessary to distinguish
  - Functions f: they give numerical transformations of parameters
    p (n) = (α,ß),   f (α,ß) = n'
    There are two functions f defined by [3] and [4]
    $f_1$ (α,ß) = (α  ß-1  9-ß  10-α), 0 < α ≤ 9, 0 < ß ≤ 9, α ≥ ß
    $f_2$ (α,0) = (α-1  9  9  10-α), 0 < α ≤ 9
  - Functions $K_i$: they yield parametric images of parameters
    p (n) = (α,ß),   p (n') = (α',ß'),  $K_i$ (α,ß) = (α',ß')

There are 13 functions $K_i$ whose domains of existence depend on the specific values α and ß as shown below

Functions $K_i$

In order to delve into the transformation process, it is necessary to know the parameters of the image. And here lies the first obstacle of Kaprekar's routine, as the image's digits are not necessarily arranged. The digits' arrangement required by the process will depend on the values of α and ß. In turn, as we have just pointed out, the structure of the transformed numbers will depend on the value of ß

a) If ß>0, n' = (α  ß-1  9-ß  10-α), p (n) = (α,ß)

The 24 permutations of the digits in a four-digit number can be grouped in 6 types based on the transposition of the extreme digits, the middle ones or both (Table 1)

Table 1. Permutations of four-digit numbers

| Type | 1 | 2 | 3 |
|---|---|---|---|
| a | $P_1$ (1 2 3 4)<br>$P_2$ (1 3 2 4)<br>$P_3$ (4 2 3 1)<br>$P_4$ (4 3 2 1) | $P_9$ (1 2 4 3)<br>$P_{10}$ (1 4 2 3)<br>$P_{11}$ (3 2 4 1)<br>$P_{12}$ (3 4 2 1) | $P_{17}$ (1 3 4 2)<br>$P_{18}$ (1 4 3 2)<br>$P_{19}$ (2 3 4 1)<br>$P_{20}$ (2 4 3 1) |
| b | $P_5$ (3 1 4 2)<br>$P_6$ (3 4 1 2)<br>$P_7$ (2 1 4 3)<br>$P_8$ (2 4 1 3) | $P_{13}$ (4 1 3 2)<br>$P_{14}$ (4 3 1 2)<br>$P_{15}$ (2 1 3 4)<br>$P_{16}$ (2 3 1 4) | $P_{21}$ (4 1 2 3)<br>$P_{22}$ (4 2 1 3)<br>$P_{23}$ (3 1 2 4)<br>$P_{24}$ (3 2 1 4) |

If we agree that $P_1$ (1 2 3 4) = (α  ß-1  9-ß  10-α), then $P_5$ = (9-ß  α  10-α  ß-1).

The most interesting thing are the possible permutations of O(n'), which come down to the framed ones above. All the arrangements where number 2 is to the left of number 1 cannot occur, because digits sorted in a descending order yields ß-1≥α, which violates the α ≥ ß restriction.

Neither can O (n') = P (3 1 2 4) = (9-ß  α  ß-1  10-α) occur, since 9-ß ≥ α, α ≥ ß-1, ß-1 ≥ 10-α → α+ß ≤ 9, α+ß ≥ 11, which is contradictory. With this there are 11 possible arrangements of O (n') left.

A function $K_i$ can be associated to each arrangement of $P_i$ after estimating α', ß'. Thus, for permutation $P_1$

α'= α-(10-α) = 2α-10,  ß'= ß-1-(9-ß) = 2ß-10

returning $K_1$ (α,ß) = (2α-10, 2ß-10)

The function's domain of existence is restricted to the values of α and ß, which lead to the corresponding permutation. In our example

α ≥ ß-1, ß-1 ≥ 9-ß, 9-ß ≥ 10-α → α ≥ ß+1,  ß ≥ 5

The 11 functions $K_i$ in $A_4$ and their domains of existence for ß≥0 are shown in Table 2.

Table 2. Functions $K_i(\alpha,\beta) = (\alpha',\beta')$ in $A_4$ with $\beta>0$

| Type | P= O(n') | $K_i(\alpha,\beta)$ | $(\alpha', \beta')$ | Existence conditions |
|---|---|---|---|---|
| 1a | $P_1$ (1 2 3 4) | $K_1$ | $(2\alpha-10, 2\beta-10)$ | $\alpha \geq \beta+1$, $\beta \geq 5$ |
|    | $P_2$ (1 3 2 4) | $K_2$ | $(2\alpha-10, 10-2\beta)$ | $\alpha \geq 6$, $\beta \leq 5$, $\alpha+\beta \geq 11$ |
| 1b | $P_5$ (3 1 4 2) | $K_5$ | $(10-2\beta, 2\alpha-10)$ | $\alpha \geq 5$, $\alpha+\beta \leq 9$ |
|    | $P_6$ (3 4 1 2) | $K_6$ | $(10-2\beta, 10-2\alpha)$ | $\alpha \leq 5$, $\alpha \geq \beta+1$ |
| 2a | $P_9$ (1 2 4 3) | $K_9$ | $[(\alpha+\beta)-9, (\alpha+\beta)]-11]$ | $\alpha \leq \beta+1$, $\alpha+\beta \geq 11$ |
|    | $P_{10}$ (1 4 2 3) | $K_{10}$ | $[(\alpha+\beta)-9, 11-(\alpha+\beta)]$ | $\alpha \geq 5$, $\beta \geq 5$, $\alpha+\beta \leq 11$ |
| 2b | $P_{13}$ (4 1 3 2) | $K_{13}$ | $[11-(\alpha+\beta), (\alpha+\beta)-9]$ | $\alpha \leq 5$, $\beta \leq 5$, $\alpha+\beta \geq 9$ |
|    | $P_{14}$ (4 3 1 2) | $K_{14}$ | $[11-(\alpha+\beta), 9-(\alpha+\beta)]$ | $\beta \leq \alpha \leq \beta+1$, $\alpha+\beta \leq 9$ |
| 3a | $P_{17}$ (1 3 4 2) | $K_{17}$ | $[(\alpha-\beta)+1, (\alpha-\beta)-1]$ | $\alpha \geq \beta+1$, $9 \leq \alpha+\beta \leq 11$ |
|    | $P_{18}$ (1 4 3 2) | $K_{18}$ | $[(\alpha-\beta)+1, 1-(\alpha-\beta)]$ | $\alpha \geq 5$, $\beta \leq 5$, $\beta \leq \alpha \leq \beta+1$ |
| 3b | $P_{21}$ (4 1 2 3) | $K_{21}$ | $[1-(\alpha-\beta), (\alpha-\beta)+1]$ | $\alpha = \beta = 5$ |

b) If $\beta = 0$, n' = $(\alpha-1\ 9\ 9\ 10-\alpha)$
There are two possible permutations of O (n')
$P_{25}$ (1 2) = (9 9 $\alpha-1$ 10-$\alpha$) y $P_{26}$ (2 1) = (9 9 10-$\alpha$ $\alpha-1$)
which lead to the parametric functions $K_{25}$ and $K_{26}$ shown in Table 3.

Table 3. Functions $K_i(\alpha,0) = (\alpha',\beta')$ in $A_4$

| P = O(n') | $K_i$ | $(\alpha',\beta')$ | Existence conditions |
|---|---|---|---|
| $P_{25}$(1 2) | $K_{25}(\alpha,0)$ | $(\alpha-1, 10-\alpha)$ | $6 \leq \alpha \leq 9$ |
| $P_{26}$ (2 1) | $K_{26}(\alpha,0)$ | $(10-\alpha, \alpha-1)$ | $0 < \alpha \leq 5$ |

## 5. Balance in $A_4$ (4-digits)

By definition, there is balance if K (n) = n = $n_E$. We will call $n_E$ constant or 1-link cycle (1-cycle).
Note that $n_E \in B_4$, since all numbers can be transformed, so $n_E$ must belong to the set which contains all the images.
In the same way, for [5]
p (n) = $(\alpha,\beta)$, p (n') = $(\alpha',\beta')$, $(\alpha',\beta') = (\alpha,\beta) = (\alpha,\beta)_E$
there is balance if, and only if, there is a given function $K_i$ that for a value $\alpha_1$, $\beta_1$ of the parameters $(\alpha,\beta)$ verifies $K_i(\alpha_1,\beta_1) = (\alpha_1,\beta_1)$, where $(\alpha_1,\beta_1)$ must belong to the domain of existence of $K_i$.

At first, no class (α,0) can be a balance class, since its image through $K_{25}$ or $K_{26}$ yields ß' ≠ 0.

The same is true for 10 of the functions $K_i$ with ß > 0. For example, for $K_{10}$ (α+ß)-9 = α, 11- (α+ß) = ß → α = -7, ß = 9 which is impossible

**But for $K_5$, associated with the permutation $P_5$ (3  1  4  2) there is in fact a solution 10-2ß = α,  2α-10 = ß → α = 6,  ß = 2 which belongs to the domain of existence of $K_5$, $K_5$ (6,2) = (6,2)  and $f_1$ (6,2) = (6  2-1  9-2  10-6) = 6174, which is Kaprekar's constant.**

Naturally, p (6  1  7  4) = (6,2) and the process is repeated ad infinitum

n: 6174 → 6174 → 6174…

(α,ß) : (6,2) → (6,2) → (6,2)…

This result is neither new nor does it demonstrate the convergence of all parametric classes to (6,2) (Prichett et al., 1981). However, the parametric transformation graph shows that all classes are iteratively transformed to (6,2).

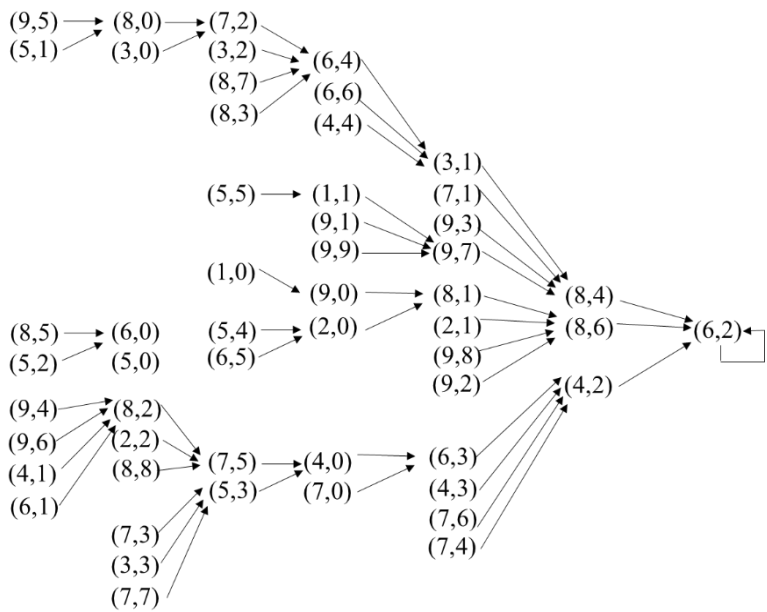

Graph 1. Tree of parametric transformations in A4

It is merely a new derivation of the constant.

## 6. Balance in $A_3$, $A_2$  and $A_5$

In the same sense, derivations should be understood in the cases of 3, 2 and 5 digits. As mere examples of the application of the $K_i$ functions to already known results.

### 6.1 $A_3$ (3-digits)

Let reference set $A_3$ consist of all natural numbers up to three digits long, excluding zero and numbers where all digits are identical. If a number has less than three digits, zeros will be added to its left.

The basic transformation equation depends on a single parameter, α:

$n = (a\ b\ c)$, $O(n) = (a_0\ b_0\ c_0)$, $0 < \alpha = a_0-c_0 \leq 9$

$f(\alpha) = (\alpha-1\ 9\ 10-\alpha) = n'$, $p(n)=\alpha$, $n'=K(n)$

For each α there is a unique image belonging to $B_3$. There are only two possible permutations of $O(n')$, $P_{12}$ and $P_{21}$, which enable two functions $K_i$

- $P_{12} = (9\ \alpha-1\ 10-\alpha)$, $6 \leq \alpha \leq 9$, $K_1(\alpha) = (\alpha-1)$
- $P_{21} = (9\ 10-\alpha\ \alpha-1)$, $0 < \alpha \leq 5$, $K_2(\alpha) = (10-\alpha)$

The two possible balance conditions $K_i(\alpha) = \alpha$, are

$K_1 : \alpha=\alpha-1$    Impossible

$K_2 : \alpha=10-\alpha \leftrightarrow \alpha=5$    Possible $\alpha_E = 5$

$K_2(5) = (5)$, $f(5) = (5-1\ 9\ 10-5) = 4\ 9\ 5$

which is the second Kaprekar's constant

There is a single tree of parametric transformations that results directly from $K_1$ y $K_2$

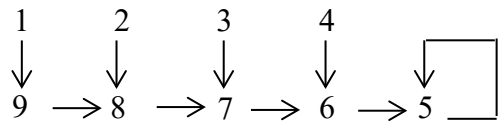

and that shows the convergence of the routine towards 495.

### 6.2 A₂ (2-digits)

Let reference set $A_2$ consist of all natural numbers up to two digits long, excluding zero and numbers with identical digits. If a number has one single digit, a zero is added to its left.

The basic transformation equation depends on a single parameter, α:

$n = (a,b)$, $O(n) = (a_0,b_0)$, $\alpha = a_0-b_0$

$n' = f(\alpha) = (\alpha-1\ 10-\alpha) \in B_2$, $p(n) = \alpha$, $n'= K(n)$

There are two possible permutations, $P_{12}$ and $P_{21}$, which enable two functions $K_i$

- $6 \leq \alpha \leq 9$, $\alpha' = 2\alpha-11$, $K_1(\alpha) = (2\alpha-11)$
- $0 \leq \alpha \leq 5$, $\alpha' = 11-2\alpha$, $K_2(\alpha) = (11-2\alpha)$

There is no value of α that verifies $K_i(\alpha_1) = (\alpha_1)$, so in set $A_2$ Kaprekar's routine cannot lead to a situation where there is balance.

### 6.3 A₅ (5-digits)

Let this reference set consist of numbers up to five digits long, excluding zero and numbers where all digits are identical.

The basic transformation equations depend on two parameters, α and ß

$n = (a\ b\ c\ d\ e)$, $O(n) = (a_0\ b_0\ c_0\ d_0\ e_0)$, $a_0 \geq b_0 \geq c_0 \geq d_0 \geq e_0$

$a_0 > e_0$, $\alpha = a_0-e_0 > 0$, $\text{ß} = b_0-d_0$, $\alpha \geq \text{ß}$

$f_1(\alpha,\text{ß}) = (\alpha\ \text{ß}-1\ 9\ 9-\text{ß}\ 10-\alpha)$  $\alpha > 0$, $\text{ß} > 0$

$f_2(\alpha,0) = (\alpha-1\ 9\ 9\ 9-\text{ß}\ 10-\alpha)$  $\alpha > 0$

For a number to be another number's image, it is necessary and sufficient to satisfy one of the following conditions:

Requirement 1) $n= (a\ b\ 9\ 8-b\ 10-a)$, $0 < a \geq b+1$

Requirement 2) $n= (a\ 9\ 9\ 9\ 9-a)$, $0 \leq a \leq 8$

The parameters of the transformed number α',ß' depend on the permutations of O (n'). The corresponding functions $K_i$ are shown in Table 4.

Table 4. Functions $K_i$ in $A_5$

| Type | P=O(n') | Existence conditions | $K_i(α,ß) = (α',ß')$ |
|---|---|---|---|
| | ß>0   $P_{1234}$ = (9  α  ß-1  9-ß  10-α),   α>ß | | |
| 1a | $P_1$ (1 2 3 4) | α≥ß+1, 5≤ß≤8 | $K_1$ (α,ß) = (α-1,α+ß-9) |
| | $P_2$ (1 3 2 4) | α≥6, 2≤ß≤5, α+ß≥11 | $K_2$ (α,ß) =[α-1,(α-ß)+1] |
| 1b | $P_5$ (3 1 4 2) | α≥5, α+ß≤9 | $K_5$ (α,ß) =[10-ß,(α-ß)-1] |
| | $P_6$ (3 4 1 2) | α≤5, α≥ß+1 | $K_6$ (α,ß) =[10-ß,9-(α+ß)] |
| 2a | $P_9$ (1 2 4 3) | α≤ß+1, α+ß≥11 | $K_9$ (α,ß) =(ß,2α-10) |
| | $P_{10}$ (1 4 2 3) | α≥5, ß≥5, α+ß≤11 | $K_{10}$ (α,ß) =[ß,(α-ß)+1] |
| 2b | $P_{13}$ (4 1 3 2) | α≤5, ß≤5, α+ß≥9 | $K_{13}$ (α,ß) =[10-ß,1-(α-ß)] |
| | $P_{14}$ (4 3 1 2) | ß≤α≤ß+1, α+ß≤9 | $K_{14}$ (α,ß) =(10-ß,10-2α) |
| 3a | $P_{17}$ (1 3 4 2) | α≥ß+1, 9≤α+ß≤11 | $K_{17}$ (α,ß) =(10-ß,2α-10) |
| | $P_{18}$ (1 4 3 2) | α≥5, ß≤5, ß≤α≤ß+1 | $K_{18}$ (α,ß) =[10-ß,(α+ß)-9] |
| 3b | $P_{21}$ (4 1 2 3) | α=ß=5 | $K_{21}$ (α,ß) =[ß,11-(α+ß)] |
| | ß=0  $P_{12}$ = (9 9 9 α-1  10-α) | | |
| 4 | $P_{25}$ (1 2) | 6<α≤9 | $K_{25}$ (α,0) =(α-1,10-α) |
| | $P_{26}$ (2 1) | 0<α≤5 | $K_{26}$ (α,0) =(10-α,α-1) |

Note that the structure in Table 4 is similar to that one corresponding to $A_4$ (Tables 2 and 3), including the functions' domains of existence. But these functions change. This is due to the fact that O(n') always starts with a 9. The functions $K_i$ act on the same parameters α and ß in both cases, but the images differ.

If the balance condition is imposed $(α_1,ß_1) = K_i (α_1,ß_1)$ there is no $K_i$ to verify it. For example, let $K_6$: α=1, ß=9, impossible since α≥ß.

Let $K_5$: α = 10-ß,  ß = (α-ß) – 1 → α = 7, ß = 3, but this solution is not valid, as it does not belong to the domain of existence (α ≥ 5, α+ß ≤ 9) of $K_5$. Actually, the image for (7,3) is (7,4).

In conclusion, in $A_5$ no pair α,ß exists which is actually its own image. Equivalently, no number exists which happens to be its own image.

### 7. **Cycles in $A_5$ and $A_2$**

There will be a cycle of r links if and only if for $(α,ß) = (α_1,ß_1)$ $K^r(α_1,ß_1) = (α_1,ß_1)$, being $K^r$ the result of operating K r times. Note that

$K^r(α,ß) = g(α,ß)$ and, in general, $g(α,ß) ≠ (α,ß)$

The condition is that there be parameters $α_1,ß_1$ such that $g(α_1,ß_1) = (α_1,ß_1)$, and that these parameters belong to the domain of existence as $g(α,ß)$.

Let us analyze sets $A_5$ and $A_2$ where there is no balance

**a)** **$A_5$**

- There exists a cycle of two links. In fact

$K^2(α,ß) = (K_{25} × K_{17})(α,ß) = K_{25}(10-ß, 2α-10) → α=5$

$K_{25}(10-ß, 0) = (9-ß, ß) = (α,ß) → α+ß=9$

$K^2(α,ß) = (α, 9-α)$, $K^2(5,4) = (5,4)$

Therefore there is a cycle of two links where one is (5,4). The other one is $K_{17}(5,4) = (6,0)$

The compatible pairs $K_i$, $K_j$, $K^2(α,ß) = (K_j × K_i)(α,ß)$ which satisfy the condition $(K_j × K_i)(α_1,ß_1) = (α_1,ß_1)$ are $K_i = K_5, K_6, K_{13}, K_{14}, K_{17}$ and $K_{18}$, $K_j = K_{25}$. All of them lead to the same result. It is also verified with $K^2(α,ß) = (K_i × K_j)(α,ß)$ for the same values. The other compatible pairs do not satisfy the conditions.

It is therefore a very simple cycle: $(6,0) \underset{K_j}{\overset{K_i}{\rightleftarrows}} (5,4)$.

Ending necessarily with the numbers

5 3 9 5 5 (6,0) ⇄ 5 9 9 9 4 (5,4)

since $f_2(6,0) = 5\ 9\ 9\ 9\ 4$ and $f_1(5,4) = 5\ 3\ 9\ 5\ 5$

- There exists no cycle with three links
- There exist two cycles of four links. In fact,

the $K^4$ satisfying the condition are derived from

$K^4(α,ß) = (K_j × K^3)(α,ß)$ with $K^3(α,ß) = (K_{17} × K_1 × K_2)(α,ß) = (19-2α+ß, 2α-14)$ valid for (8,4) and (8,3)

$K_j = K_5 → K^4(α,ß) = (24-2α, 32-4α+ß)$   α = 8

$24-2α = α → α = 8$, $K^4(8,ß) = (8,ß)$

As $K^3$ can only take two arguments (8,4) and (8,3), there are two solutions which lead to two cycles of four links:

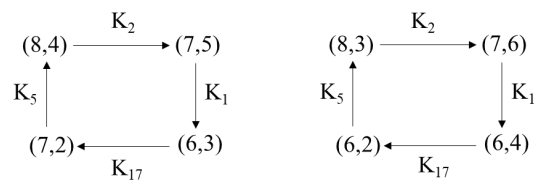

Needless to say, any of the four $K_4$ formed by rotation of operators $K_i$ leads to the same cycle. In terms of numbers

$$\begin{array}{l}\rightarrow 7\,1\,9\,7\,3\,(8,4) \rightarrow 8\,3\,9\,5\,2\,(7,5) \\ \phantom{\rightarrow} 6\,2\,9\,6\,4\,(7,2) \leftarrow 7\,4\,9\,4\,3\,(6,3) \leftarrow\end{array} \qquad \begin{array}{l}\rightarrow 6\,1\,9\,7\,4\,(8,3) \rightarrow 8\,2\,9\,6\,2\,(7,6) \\ \phantom{\rightarrow} 6\,3\,9\,5\,4\,(6,2) \leftarrow 7\,5\,9\,3\,3\,(7,4) \leftarrow\end{array}$$

Prichett (1978) showed the existence of these cycles but not their algebraic derivation.

- There is no other combination of $K_i$ in $A_5$ that can generate cycles different from the ones stated above.

**b)** **$A_2$**

No $K^2$ or $K^3$ satisfies the conditions $K^r(\alpha_1) = \alpha_1$

For instance,

$K^2(\alpha) = (K_1 \times K_1)(\alpha) = (4\alpha-33)$ implies $\alpha = 11$
$K^2(\alpha) = (K_1 \times K_2)(\alpha) = (11-4\alpha)$, $\alpha = 11/5$
$K^2(\alpha) = (K_2 \times K_2)(\alpha) = (4\alpha-11)$, $\alpha = 11/3$
$K^3(\alpha) = K_2^3(\alpha) = (33-8\alpha)$, $\alpha = 11/3$

impossible in all cases, therefore there aren't any cycles of two, three, or four links.

However,

$K^5(\alpha) = K_1^2 \times K_2^3(\alpha) = (99-32\alpha) \rightarrow \alpha = 3$

There is consequently a cycle of five links

$$\begin{array}{ccc} & K_2 & \\ (5) & \rightarrow & (1) \\ K_2 \nearrow & & \searrow K_2 \\ (3) & \leftarrow (7) \leftarrow (9) \\ & K_1 \quad K_1 & \end{array}$$

which articulate the rest of the parameters (2), (4), (6) and (8) according to $K_1$ and $K_2$. Any 2-digit number will become one of the five following numbers

$$\begin{array}{ccc} 27\,(5) & \rightarrow & 45\,(1) \\ \nearrow & & \searrow \\ 63\,(3) \leftarrow & 81\,(7) \leftarrow & 09\,(9) \end{array}$$

which are shown alongside their parameter.

**8.** **The numeral system**

Obvious as it may seem, it is worth noting the importance of the numeral system. This aspect has been and is being intensively investigated by Walden (2005), Hanover (2017), Yamagami (2018), Devlin & Zeng (2020), Wang & Lu (2021),… Here we simply present an example of the application of the functions $K_i$.

Let us assume we are using a numeral system in base 100. Then, a 4-digit number in base 10 would become a 2-digit number in base 100 (not including decimal numbers from 9901 to 9999, the latter being discarded for having all identical digits).

If we represent number 99 by z, the result is

$n = (a\ b)_{100}$, $0 \leq a \leq z$, $0 \leq b \leq z$
$O(n) = (a_0\ b_0)_{100}$, $\alpha = a_0-b_0$, $0 < \alpha \leq z$

And the basic functions written in base 100 are identical to those of $A_2$ in base 10

$f(\alpha) = (\alpha-1\ \ 10-\alpha)_{100}$, $n' = \dot{z}$
$K_1(\alpha) = (2\alpha-11)_{100}$, $2\alpha > 11$

$K_2(\alpha) = (11-2\alpha)_{100}$, $2\alpha < 11$

thus making balance impossible

$K_1(\alpha) = \alpha \rightarrow \alpha = 2\alpha-11 \rightarrow \alpha = 11 > z$
$K_2(\alpha) = \alpha \rightarrow \alpha = 11-2\alpha \rightarrow \alpha = 11/3$
The system leads to a cycle consisting of 50 links

### 9. Application of the model to other transformation processes

Let the following process start with any number n with no repeated digits. As in Kaprekar's process, n's digits are sorted in descending order, yielding another number X. Another number Y is obtained by translocation the extreme digits in X, then n' = X-Y is calculated. The process is iterated repeatedly. For example, if n = 8 0 7 2, X= 8 7 2 0, Y = 0 7 2 8, n' = 7 9 9 2, n'' = 6 9 9 3, n''' = 5 9 9 4, $n^{4)}$ = 4 9 9 5

For any n∈$A_w$ the basic function will be
$O(n) = (a_0\ b_0\ \ldots z_0)$, $a_0 \geq b_0 \geq \ldots \geq z_0$, $\alpha = a_0 - z_0 > 0$

$f(\alpha) = (\alpha-1\ \overset{s}{9\ldots9}\ 10-\alpha)$, $s = w-2 \geq 0$

For w=2 and 3 this process matches Kaprekar's. The functions $K_i$ in $A_w$, w > 2 are
$K_1(\alpha) = (\alpha-1)$, $6 \leq \alpha \leq 9$
$K_2(\alpha) = (10-\alpha)$, $0 < \alpha \leq 5$
valid regardless of the number of digits greater than 2 that n has.
There is balance when

$K_2(\alpha) = \alpha \rightarrow \alpha = 5$ 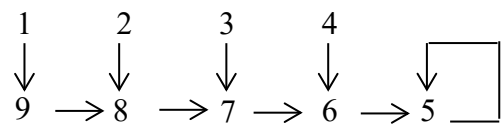

which is a generalization of the second Kaprekar's constant for numbers with any number of digits greater than 2.

The tree of parametric transformations is identical to that of the 3-digit Kaprekar, as a consequence of the fact that the $K_i$ functions are identical

```
1    2    3    4
↓    ↓    ↓    ↓    ↓
9 → 8 → 7 → 6 → 5 ⌐
```

this shows the convergence of the tree.

Perhaps the $K_i$ functions can be useful to study other numerical transformation processes.

### 10. Discussion

We start with set $A_4$, consisting of all 4-digit natural numbers, excluding those with repeated digits. What do all the numbers yielding the same image have in common? Two parameters α and ß, which we define after sorting the number's digits. That is, in terms of transformation, any number is marked by these two parameters. All numbers having the same parameters will yield the same image.

The basic functions $f_1$ and $f_2$ establish the numerical image of any number of given parameters:

$f_1(\alpha,\beta) = (\alpha \; \beta-1 \; 9-\beta \; 10-\alpha)$, $\alpha \geq \beta$, $\beta > 0$

$f_2(\alpha,\beta) = (\alpha-1 \; 9 \; 9 \; 10-\alpha)$, $\alpha \geq \beta$, $\beta = 0$

These functions have already been established by Prichett et al. (1981).

These functions impose severe restrictions to the images. All of them are multiples of nine. What's more, the structure of their digits is subject to the parameter ß of the number to be transformed being greater than or equal to zero. If ß > 0, adding the first and fourth digits yields 10, and adding the second and third ones yields 8. If ß = 0, adding the extreme digits returns 9, as happens with the middle ones. But also for a number to belong to $B_4$ it must have an anti-image. The transformation is an epimorphism of $A_4$ in $B_4$, the latter being a subset of those multiples of 9 that satisfy the aforementioned requirements and that have anti-image. The necessary and sufficient conditions for a number to belong to $B_4$ are a novel contribution. This aspect will be key in the analysis of the transformation trees architecture.

If we want to make progress in the transformation process, it is mandatory to know the parameters of the transformed number. And here lies the main obstacle in Kaprekar's routine, as the image's digits are not necessarily arranged. The digits' arrangement required by the transformation will depend on the values of α and ß. Thus, if $\alpha \geq 5$, ß > 0 and $\alpha+\beta \leq 9$, then the descending arrangement of transformed number n' is (9-ß α 10-α ß-1) and its parameters are α' = 10-2ß, ß' = 2α-10.

Any 4-digit number can have up to 24 permutations, but because it is $\alpha \geq \beta$, with ß > 0, only 11 of them are possible in the ordered image. With ß = 0 there are two more possible permutations.

Each possible permutation can have a function $K_i$ attached whose argument (α,ß) is the parameters of the number to be transformed, and whose image (α',ß') is the parameters of the transformed number. Therefore, in our previous example

$K_5(\alpha,\beta) = (10-2\beta, 2\alpha-10)$, $\beta > 0$, $\alpha \geq 5$, $\alpha+\beta \leq 9$

With these functions $K_i$ we are able to turn sets of numbers with the same image into the set of transformed numbers. This way, we simplify the transformation process by going from 9990 numbers to 54 classes.

Furthermore, with $K_i$ functions we have a key tool to approach the problem of equilibrium and the nature of cycles from a new angle and the algebraic architecture of the transformations, our final objective.

By definition, there is balance when the image of a number is that number itself, i.e., K(n) = n. Other authors identify equilibrium with the existence of a 1-link cycle (1-cycle). It is merely a question of nomenclature. Equivalently, it is a necessary and sufficient balance condition that there be a function $K_i$ such that for a specific value $\alpha_1,\beta_1$ returns the parameters $K_i(\alpha_1,\beta_1) = (\alpha_1,\beta_1)$,

In set $A_4$ of 4-digit numbers only one function satisfies this condition – the aforementioned function $K_5$.

$\alpha = 10-2\beta$, $\beta = 2\alpha-10 \rightarrow \alpha = 6$, $\beta = 2$

$K_5(6,2) = (6,2)$ and applying the function $f_1$

$f_1(6,2) = 6174$

which is the renowned Kaprekar's constant.

The developed methodology allows to approach the balance in the set of numbers $A_w$ with w digits. In $A_3$ the basic transformation equation depends on a single parameter $f(\alpha) = (\alpha-1 \quad 9 \quad 10-\alpha)$ and two functions $K_i$ exist, $K_1(\alpha) = (\alpha-1)$ and the function $K_2(\alpha) = (10-\alpha)$ that determines balance, since $K_2(5) = (5)$, $f(5) = (4 \quad 9 \quad 1)$, resulting in the second Kaprekar constant. Convergence is obvious from these functions.

In $A_2$ the process depends also on a single parameter with $f(\alpha) = (\alpha-1 \quad 10-\alpha)$ and two functions $K_i$, none of which satisfies the balance condition.

In $A_5$ the algorithm depends again on two parameters with two basic functions $f_1$ y $f_2$ and 13 functions $K_i$, none of which satisfies the balance condition.

Prichett et al. (1981) establishes that for $w = \dot{2}$, $w \geq 6$ there are at least two different constants (1-cycle). The same occurs for $w = \dot{2} + 1$, $w \geq 15$. The work also establishes characteristics of these constants, since there are multiple solutions. Also these authors (Lapenta et al., 1979) use algorithms to determine some constants.

$K_i$ functions offer an alternative way to determine exactly these constants, although this is only operative if the number of $K_i$ functions is not very large.

The developed tools $K_i$ enable also the study of cycles of more than 1-link. The existence of cycles entails the fact that the repetition of the transformation process provides parameters identical to the original ones. Officially, there is a cycle of r links if and only if for $(\alpha,\beta) = (\alpha_1,\beta_1)$ $K^r(\alpha_1,\beta_1) = (\alpha_1,\beta_1)$, being $K^r$ the result of applying the transformation algorithm r times. As we have already stated ultimately, balance would make the specific case $r = 1$.

In the cases studied we have been able to determine $K^r$ functions that determine the r-cycles, $r \geq 2$. Thus, in $A_5$, $K^2(\alpha,\beta) = (\alpha,9-\alpha)$ valid for $\alpha+\beta = 9$ y $\alpha = 5$, giving rise to the cycle $(5,4) \rightleftarrows (6,0)$. Besides $K^4(\alpha,\beta) = (24-2\alpha, 32-4\alpha+\beta)$ admits two solutions $(8,4)$ y $(8,3)$ giving rise to the two 4-ciclos specified in section 7a. In $A_2$ the function $K^5(\alpha) = (99-32\alpha)$, determinates a 5-cycle.

The nature of the r-cycles is independent of the number of parameters. It depends on the structure of the $K_i$ and this of the basic functions f.

It would seem obvious to state that the structure of the transformation tree generated by Kaprekar's process is influenced by the numeral system. Both the number of parameters and the basic functions depend on it. It is an important topic that is being intensively studied as indicated in section 8. Let us take as an example of application of functions $K_i$ the case $A_4$, which when represented in base 100 becomes a case $A_2$. In this base, the system has one single parameter, as well as identical functions to those of $A_2$ in base 10, albeit different domains of existence, presenting a cycle of 50 links.

We hope that the functions $K_i$ will be useful for further study of other numerical transformation processes. As an example, we introduce a generalization of the second Kaprekar's constant for numbers of w digits, $w > 2$.

The number of transformation trees in each case will depend on the coupling of functions $K_i$, that is, the existence of functions whose argument matches the image of another one. Both this topic and the architecture of trees will be dealt with in the second part of this paper (Nuez, 2021).

## 11. <u>References</u>